\documentclass[11pt,a4paper,reqno]{amsart}
\usepackage{amsfonts,amssymb,amscd,amsmath,enumerate,verbatim,calc}

\newtheorem{theorem}{Theorem}

\begin{document}

\title[Brouwer's fixed point theorem and Sperner's lemma]
{On the (in)equivalence of Brouwer's fixed point theorem and Sperner's lemma}

\author[J.~Minagawa]{Junichi Minagawa}

\address{%
Faculty of Economics,
Chuo University,
742-1 Higashinakano,
Hachioji,
Tokyo 192-0393,
Japan}

\email{minagawa@tamacc.chuo-u.ac.jp}

\begin{abstract}
  We consider Brouwer's fixed point theorem and Sperner's lemma in one dimension.
  We present a proof of the Brouwer theorem using the Sperner lemma, and vice versa.
  However, we also show that they are not equivalent, because the Sperner lemma holds
  in the ordered field of rational numbers, whereas proving the Brouwer theorem requires
  the property of completeness.
\end{abstract}

%\subjclass[2010]{26A03; 05C15}


\maketitle

\noindent
Brouwer's fixed point theorem
\cite{Brouwer}
%(1912)
and Sperner's lemma
\cite{Sperner}
%(1928)
have had profound implications for
various fields.
%game theory and economics.
It is well established that Brouwer's fixed point theorem can be proved by using Sperner's lemma
(see \cite{Knaster}, \cite[p.~28]{Border}),
%(see Knaster et al.~1929; Border 1985, p.~28),
and conversely, Sperner's lemma can also be proved
by using Brouwer's fixed point theorem (see \cite{Yoseloff}, \cite{Park}).
%(see Yoseloff 1974; Park 2003).
Thus, these two theorems are often considered equivalent (see, e.g.,
\cite{Yoseloff}, \cite{Garcia}, \cite{Park}, \cite{Voorneveld}, \cite{Idzik}).
%(see, e.g., Yoseloff 1974; Garcia and Zangwill 1981; Park 2003; Voorneveld 2017;
%Idzik et al.~2021).
However, the existing proofs of Brouwer's fixed point theorem via Sperner's lemma make use of other
results as well, such as the Bolzano--Weierstra{\ss} theorem.
This suggests that Brouwer's fixed point theorem cannot be derived from Sperner's lemma alone,
implying that they are not equivalent.%
\footnote{%
  From a purely logical point of view, Brouwer's fixed point theorem and Sperner's lemma
  are not equivalent in general.
  It is known in the context of reverse mathematics that Sperner's lemma is
  provable in an axiom system called $\mathrm{RCA}_{0}$
  (see \cite[p.~149]{Simpson}),
%  (see Simpson 2009, p.~149),
  whereas Brouwer's fixed point theorem for the unit interval is provable in $\mathrm{RCA}_{0}$ but
  for the unit square it is not provable in $\mathrm{RCA}_{0}$
  (see \cite{Shioji}).
%  (see Shioji and Tanaka 1990).
}

The aim of this note is to show that this is the case, using elementary methods.
In the following discussion, we focus on the one-dimensional versions of Brouwer's fixed point
theorem and Sperner's lemma.

\begin{theorem}[Brouwer's fixed point theorem in one dimension]
Suppose that $f$ is a continuous function from a closed interval $[a, b] \subset \mathbb{R}$
to itself.
Then there exists a fixed point $x \in [a, b]$ such that $f (x) = x$.
\end{theorem}

\begin{theorem}[Sperner's lemma in one dimension]
Let $V = \{ v_{0}, v_{1}, \dots, v_{n} \} \subset \mathbb{R}$, where $v_{0} < v_{1} < \dots < v_{n}$
and $n \geq 1$.
Let $\kappa : V \to \{ 0, 1 \}$ such that $\kappa (v_{0}) = 0$ and $\kappa (v_{n}) = 1$.
Then there exists an $i \in \{ 1, 2, \dots, n\}$ such that $\kappa (\{ v_{i-1}, v_{i} \}) =
\{ 0, 1 \}$.
\end{theorem}

First, we give a proof of the Brouwer theorem using the Sperner lemma, and then give a proof of
the Sperner lemma using the Brouwer theorem.
The proofs presented here are not merely one-dimensional adaptations of existing ones,
as cited above.
Instead, they draw inspiration from the work of Gale
\cite{Gale}
%(1979)
on the game of Hex
(for expositions of this work, see, e.g., \cite{Vohra}, \cite{Karlin}).
%(for expositions of this work, see, e.g., Vohra 2005; Karlin and Peres 2017).

\begin{proof}[Proof of the Brouwer theorem using the Sperner lemma]
Suppose that a continuous function $f: [a, b] \to [a, b]$ has no fixed point.
Define a function $g (x) = f (x) - x$; then, the function $g$ is continuous on $[a, b]$ and
is never zero.
Thus, by the extreme value theorem, there is an $\varepsilon > 0$ such that $| g (x) | \geq
\varepsilon$ for all $x \in [a, b]$.
Moreover, any continuous function on a closed and bounded interval is uniformly continuous,
so there exists a $\delta > 0$ such that $| x - x' | < \delta$ implies
$| g (x) - g (x') | < \varepsilon$.
Without loss of generality, we can take $\delta < \varepsilon$.
Additionally, by the Archimedean property, we can choose a positive integer $n$ such that
$(b - a) < n \delta$.

Suppose that $a = v_{0} < v_{1} < \dots < v_{n} = b$ where $v_{i} = a + i (b - a) / n$.
If $f (v_{j}) - v_{j} \geq \varepsilon$, then we set $\kappa (v_{j}) = 0$.
If $\kappa (v_{j}) \neq 0$, then $| g (v_{j}) | \geq \varepsilon$ leads to
$f (v_{j}) - v_{j} \leq - \varepsilon$ and in this case, we set $\kappa (v_{j}) = 1$.
Since $f (v_{j}) \in [a, b]$, we must have $f (v_{0}) - v_{0} \geq \varepsilon$
and $f (v_{n}) - v_{n} \leq - \varepsilon$; thus, $\kappa (v_{0}) = 0$ and $\kappa (v_{n})
= 1$.
Then, by the Sperner lemma, there exists an $\ell \in \{ 1, 2, \dots, n\}$ such that
$\kappa (\{ v_{\ell-1}, v_{\ell} \}) = \{ 0, 1 \}$,
specifically where $\kappa (v_{\ell-1}) = 0$ and $\kappa (v_{\ell}) = 1$; thus,
$f (v_{\ell-1}) - v_{\ell-1} \geq \varepsilon$ and $f (v_{\ell}) - v_{\ell} \leq - \varepsilon$.
Therefore, we obtain $f (v_{\ell-1}) - f (v_{\ell}) \geq 2 \varepsilon - (v_{\ell} - v_{\ell-1}) =
2 \varepsilon - (b - a) / n > 2 \varepsilon - \delta > \varepsilon$,
which gives a contradiction, since $| v_{\ell-1} - v_{\ell} | = (b - a) / n < \delta$
implies $| g (v_{\ell-1}) - g (v_{\ell}) | = | f (v_{\ell-1}) - f (v_{\ell}) + v_{\ell} - v_{\ell-1} | <
\varepsilon$.
\end{proof}

\begin{proof}[Proof of the Sperner lemma using the Brouwer theorem]
Suppose that there exists no $i \in \{ 1, 2, \dots, n\}$ such that $\kappa (\{ v_{i-1}, v_{i} \}) =
\{ 0, 1 \}$.
Let $A = \{ v_{j} \in V | \kappa (v_{j}) = 0 \}$, and let $B = \{ v_{j} \in V |
\kappa (v_{j}) = 1 \}$.
Define a function $f$ from $V$ to itself as follows:
\begin{equation}
  f (v_{j}) =
  \begin{cases}
    v_{j+1} \; \; \text{if} \; v_{j} \in A  \\
    v_{j-1} \; \; \text{if} \; v_{j} \in B,
  \end{cases}
\end{equation}
where we assume that $v_{-1} = c < v_{0}$ and $v_{n+1} = d > v_{n}$.
We now verify that $f (v_{j}) \in V$.
If $v_{j} \in A$ and $v_{j+1} \notin V$, then we must have $v_{j} = v_{n}$;
hence, $v_{j} \in B$, which contradicts $v_{j} \in A$.
If $v_{j} \in B$ and $v_{j-1} \notin V$, a contradiction arises in a similar way.

We extend $f$ to $\hat{f} : [v_{0}, v_{n}] \to [v_{0}, v_{n}]$ as follows.
For any $x \in [v_{0}, v_{n}]$, we can express $x$ as a convex combination of some
$v_{k - 1}, v_{k} \in V$ where $k \in \{ 1, 2, \dots, n\}$, that is,
$x = \lambda v_{k - 1} + (1 - \lambda) v_{k}$ and $0 \leq \lambda \leq 1$.
Define $\hat{f} (x) = \lambda f (v_{k - 1}) + (1 - \lambda) f (v_{k})$;
then, the function $\hat{f}$ is continuous.
Hence, by the Brouwer theorem, there exists a fixed point $x \in [v_{0}, v_{n}]$ such that
$\hat{f} (x) = x$.
When $x$ is a fixed point, it follows that $\lambda v_{k - 1} + (1 - \lambda) v_{k} =
\lambda f (v_{k - 1}) + (1 - \lambda) f (v_{k}) = \lambda v_{k - 1 + e} +
(1 - \lambda) v_{k + e'}$, where $e, e' \in \{ 1, - 1 \}$.
Suppose that $v_{k - 1}, v_{k} \in A$; then $e = e' = 1$.
Thus, we have $\lambda v_{k - 1} + (1 - \lambda) v_{k} <
\lambda v_{k - 1 + e} + (1 - \lambda) v_{k + e'}$, a contradiction.
Similarly, we obtain a contradiction for the case $v_{k - 1}, v_{k} \in B$.
\end{proof}

Next, we establish that the Sperner lemma does not imply the Brouwer theorem,
for the following reason.
When we replace $\mathbb{R}$ with $\mathbb{Q}$, the former continues to hold but the latter fails.
Therefore, the Sperner lemma is not equivalent to the Brouwer theorem, just as the Archimedean
property is not equivalent to the Dedekind completeness property
because the ordered field of rational numbers satisfies the Archimedean property but does not
satisfy the Dedekind completeness property
(see, e.g., \cite{Propp}).
%(see, e.g., Propp 2013).

\begin{proof}[Proof that the Sperner lemma does not imply the Brouwer theorem]
Suppose that the underlying ordered field is taken to be the rationals $\mathbb{Q}$.
Then, it is clear that the Sperner lemma is still true:
Let $V = \{ v_{0}, v_{1}, \dots, v_{n} \} \subset \mathbb{Q}$, where $v_{0} < v_{1} < \dots < v_{n}$
and $n \geq 1$.
Since $\kappa (v_{0}) = 0$ and $\kappa (v_{n}) = 1$,
there must be an $i \in \{ 1, 2, \dots, n\}$ such that
$\kappa (v_{i-1}) = 0$ and $\kappa (v_{i}) = 1$.
However, the Brouwer theorem is false.
An example is the function $f$ defined on $[1, 2] \cap \mathbb{Q}$ by
\begin{equation}
  f (x) =
  \begin{cases}
    2 \; \; \text{if} \; 1 \leq x < \sqrt{2}  \\
    1 \; \; \text{if} \; \sqrt{2} < x \leq 2,
  \end{cases}
\end{equation}
which is continuous but has no fixed point.
\end{proof}

The proof of the Sperner lemma presented earlier relies solely on the Brouwer theorem.
On the other hand, the proof of the Brouwer theorem not only uses the Sperner lemma but also
other results.
Proving the Brouwer theorem, in fact, requires the property of completeness.
In this regard, the following statements are known to be equivalent
(see, e.g., \cite{Riemenschneider}, \cite{Propp}, \cite{Ovchinnikov}):
%(see, \eg, Riemenschneider 2001; Propp 2013; Ovchinnikov 2021):
\begin{enumerate}
\item[(i)] Dedekind completeness property.
\item[(ii)] Brouwer's fixed point theorem in one dimension.
\item[(iii)] Extreme value theorem.
\item[(iv)] Archimedean property $+$ Any continuous function on a closed and bounded interval is
  uniformly continuous.
\item[(v)] Bolzano--Weierstra{\ss} theorem in one dimension.
\end{enumerate}
Here, equivalence means that any one of these statements can be taken as an axiom for
completeness of the real number system.
Notably, we invoked (iii) and (iv) in the proof of the Brouwer theorem using the Sperner lemma.
Finally, we refer to
\cite[p.~119]{Vohra}
%Vohra (2005, p.~119)
for a more standard proof, which uses (v) and
the Sperner lemma to prove the Brouwer theorem.

\end{document}